\documentclass[12pt,a4paper,twoside,final,notitlepage, leqno]{article}
\usepackage[english]{babel}
\usepackage[T1]{fontenc}
\usepackage{epsfig, graphicx, amssymb}
\usepackage{amsmath,amsthm,epsfig,amsfonts}
\usepackage{float}
\usepackage{color}
\def\br {\break}

\newcommand{\moneq}{\vspace*{-7pt} \begin{equation} \displaystyle }
\newcommand{\moneqstar}{\vspace*{-6pt} \begin{equation*} \displaystyle }
\newcommand{\monendstar}{\vspace*{-6pt} \end{equation*}   }
\newcommand{\monend}{\vspace*{-7pt} \end{equation}   }
\newcommand{\moneqarraystar}{ \begin{eqnarray*} \displaystyle }
\newcommand{\monendarraystar}{ \end{eqnarray*}   }

\newcommand{\dd}{{\rm d}}

\definecolor{vertfonce}{rgb}{0.0, 0.5, 0.0}

\def\section*#1{}

\usepackage{fancyhdr}
\renewcommand{\headrulewidth}{0pt}

\begin{document}

\fancypagestyle{plain}{ \fancyfoot{} \renewcommand{\footrulewidth}{0pt}}
\fancypagestyle{plain}{ \fancyhead{} \renewcommand{\headrulewidth}{0pt}}

~

  \vskip 2.1 cm

\centerline {\bf \LARGE Simpson's quadrature for a}

\bigskip

\centerline {\bf \LARGE nonlinear variational  symplectic scheme}

 \bigskip  \bigskip \bigskip

\centerline { \large    Fran\c{c}ois Dubois$^{ab}$ and Juan Antonio Rojas-Quintero$^{c}$ }

\smallskip  \bigskip

\centerline { \it  \small
  $^a$   Laboratoire de Math\'ematiques d'Orsay, Facult\'e des Sciences d'Orsay,}

\centerline { \it  \small   Universit\'e Paris-Saclay, France.}

\centerline { \it  \small
$^b$    Conservatoire National des Arts et M\'etiers, LMSSC laboratory,  Paris, France.}

\centerline { \it  \small CONAHCYT $-$ Tecnol\'ogico Nacional de M\'exico/I.T. Ensenada, Ensenada 22780, BC,
Mexico}


\bigskip  \bigskip

\centerline {20 October  2023
{\footnote {\rm  \small $\,$ This contribution is published in the
Springer Proceedings in Mathematics and Statistics, volume~433,
{\it Finite Volumes for Complex Applications X, volume 2,
Hyperbolic and Related Problems}, E. Franck {\it et al.} (eds.), 
pages 83-92, 2023. 
It has been  presented to the
10th international conference on Finite Volumes for Complex Applications, 
Strasbourg, 30 october - 03 november 2023.}}}




\bigskip  \bigskip
\noindent {\bf \large Abstract}

\noindent
We propose a variational  symplectic numerical method for the time integration of
dynamical systems issued from the least action principle.
We  assume a quadratic internal interpolation of the state between two time steps 
and we approximate the action in one  time step by the Simpson's quadrature formula.
The resulting scheme is nonlinear and symplectic. 
First numerical experiments concern a nonlinear pendulum and 
we have observed experimentally very good convergence properties.

\bigskip \bigskip    \noindent {\bf \large    1) \quad  Introduction} 

\smallskip \noindent
The principle of least action is a key point for establishing evolution equations
and  partial differential equations, from classical to quantum mechanics and electromagnetisms
\cite{Ar74,So70}.
An important application of this principle is proposed
with the finite element method 
and it is used for  engineering applications 
since the 1950's.
The  principle of least action is also the starting point  for the conception of
symplectic numerical schemes for dynamical systems (see {\it e.g.} \cite{HWL06,KMOW00,SS92}). 
In particular, the Newmark scheme \cite{Ne59}, which is very popular in structural dynamics, is symplectic. 
In \cite{Du23}, we have proposed a linear Simpson sympletic scheme that extends the previous
works \cite{KMOW00,SS92} for quadratic interpolation. 

\smallskip \noindent
In this contribution, we study a nonlinear Simpson symplectic scheme, as an alternative to Newmark's scheme.
One motivation is to solve stiff problems in robotics \cite{RVS21}.
We first recall  in Section 2 the fundamental statements relative to Lagrangian
and Hamiltonian mechanics and focus our attention on  the nonlinear pendulum.
Then the classical discrete dynamical system  obtained directly with a second order discretization
of the Lagrangian is  presented  in  Section 3.
Then we define and study in Section 4 a  symplectic Simpson numerical scheme based on
a quadrature with internal quadratic interpolation. 
First numerical results are presented in Section 5 
and the exact solution of the nonlinear pendulum is recalled in the Annex. 

\bigskip    \noindent {\bf \large    2) \quad  Continuous approach} 

\smallskip \noindent
We consider a dynamical system described by a state $ \, q(t) \, $
of constant mass $ \, m \, $ composed by a simple real variable 
to fix the ideas,  and for $ \, 0 \leq t \leq T $.
The continuous action~$ \, S_c \, $ introduces a Lagrangian~$ \, L \, $
\moneq \label{Lagrangien-continu} 
L \Big( q, \,  {{\dd q}\over{\dd t}} \Big)  = {m\over2} \, \Big( {{\dd q}\over{\dd t}} \Big)^2 - V(q) 
\monend
and we have

\vskip -.8 cm 
\moneq \label{action-continue} 
S_c = \int_0^T L \Big( {{\dd q}\over{\dd t}} ,\, q(t) \Big) \, \dd t . 
\monend  
The trajectories associated with the extremals of the action satisfy the Euler-Lagrange equations
%
$ \, {{\dd}\over{\dd t}} \big( {{\partial L}\over{\partial (  {{\dd q}\over{\dd t}} ) }} \big) =  {{\partial L}\over{\partial q}} $. 
%
With the Lagrangian proposed in (\ref{Lagrangien-continu}), the differential equation
\moneq   \label{oscillateur-nonlineaire}
{{\dd}\over{\dd t}} \Big( m \, {{\dd q}\over{\dd t}} \Big) +  {{\partial V}\over{\partial q}} = 0  
\monend
of Newtonian mechanics is recovered. 
With the momentum $ \, \smash{ p \equiv  {{\partial L}\over{\partial (  {{\dd q}\over{\dd t}} ) }} } = m \,  {{\dd q}\over{\dd t}} \, $
and the Hamilton function $ \, H(p,\, q) \equiv p \, {{\dd q}\over{\dd t}}  - L \big( q, \,  {{\dd q}\over{\dd t}} \big) $, 
we obtain the  first order system of Hamilton's equations
$ \, {{\dd p}\over{\dd t}} +  {{\partial H}\over{\partial q}} = 0 $,
$ \, {{\dd q}\over{\dd t}} -  {{\partial H}\over{\partial p}} = 0 $.
In the case of the Lagrangian function introduced in~(\ref{Lagrangien-continu}),
we have $ \,  H(p,\, q) = {1\over{2\, m}} \, p^2 + V(q) \, $ and 
\moneq   \label{hamilton-continu}
{{\dd p}\over{\dd t}} +  {{\partial V}\over{\partial q}} = 0 \,,\,\,  {{\dd q}\over{\dd t}} - {1\over{m}} \, p = 0 . 
\monend
%
In this contribution, we consider the case of  the nonlinear pendulum that corresponds to the  potential
$ \, \smash{ V(q) = m \, \omega^2 \, (1 - \cos q ) } $. Then we have
$ \, {{\dd p}\over{\dd t}} +  m \, \omega^2 \, \sin q = 0 \, $ and finally the second order dynamics 
$ \, {{\dd^2 q}\over{\dd t^2}} +  \omega^2 \, \sin q = 0 $.
An analytical solution of this problem is established in Annex. This  exact solution
is also shown  in Figure 1. 
%

\fancyhead[EC]{\sc{Fran\c{c}ois Dubois and Juan Antonio Rojas-Quintero}}
\fancyhead[OC]{\sc{Simpson's quadrature for a nonlinear variational  symplectic scheme}}
\fancyfoot[C]{\oldstylenums{\thepage}}

\centerline {\includegraphics[width=1.0 \textwidth]{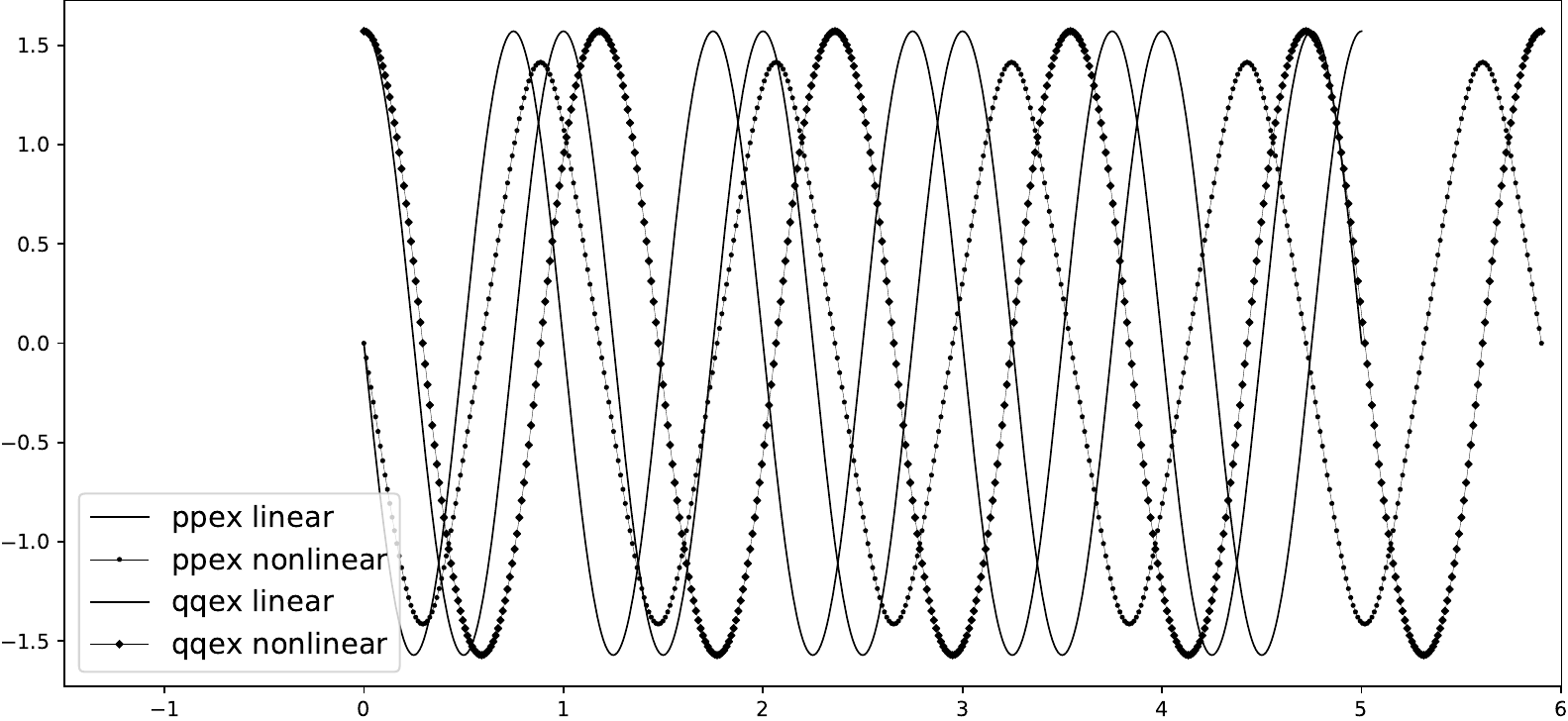}}

\noindent
Figure 1. Typical evolution for five periods of an nonlinear pendulum satisfying the dynamics  (\ref{hamilton-continu})
with $ \, \smash{ V(q) = m \, \omega^2 \, (1 - \cos q ) } $, $\, \omega = 2 \, \pi $, 
$ \, p(0) = 0 \, $ and $ \, q(0) = {{\pi}\over{2}} $. 
We observe that the nonlinear period is not given by the linear evaluation $\, T = {{2 \pi}\over{\omega}} \, $ anymore
but follows the relation~(\ref{periode}). 

\bigskip    \noindent {\bf \large    3) \quad  Newmark scheme} 

\smallskip \noindent
In this section we recall the classical Newmark scheme \cite{Ne59} which we use as the reference scheme for benchmarking purposes.
A discretization of the relation (\ref{action-continue}) is obtained  by splitting the interval
$ \, [0 ,\, T ] \, $  into $ \, N \, $ elements and we set $\, h = {{T}\over{N}} $.
At the discrete time $ \, t_j = j \, h $, an approximation  $\, q_j \, $ of
$ \, q(t_j) \, $ is introduced and a discrete form $ \, S_d \, $  of the continuous action $ \, S_c \, $ 
can be defined: 
$ \,  S_d = \sum_{j=1}^{N-1} L_N (q_j ,\, q_{j+1} ) $.
%
The discrete Lagrangian  $ \, L_N (q_\ell ,\, q_r ) \, $
is derived from the  relation  (\ref{Lagrangien-continu})  
with a centered  finite difference approximation 
$ \,  {{\dd q}\over{\dd t}} \simeq {{q_r - q_\ell}\over{h}} $ 
and a midpoint  quadrature formula
$ \, \int_0^h V \big( q(t) \big) \, \dd t \simeq h \, V \big( {{q_\ell + q_r}\over2} \big) $:
\moneq \label{lagrangien-discret} 
L_N (q_\ell ,\, q_r ) = {{m\,h}\over2} \, \Big( {{q_r - q_\ell}\over{h}} \Big)^2 - h \, V \Big( {{q_\ell + q_r}\over2} \Big) .
\monend
We observe that  $ \, S_d = \cdots  + \,  L_N (q_{j-1} ,\, q_j ) +  L_N (q_j ,\, q_{j+1} )  \, + \cdots $.
Then the discrete Euler Lagrange equation is obtained when the action is stationary,
$ \, \delta S_d = 0 \, $ for  an arbitrary variation $ \, \delta q_j \, $
of the discrete variable $ \, q_j $. It takes the form 
\moneq \label{Euler-Lagrange-discret} 
{{\partial L_N}\over{\partial q_r}}  (q_{j-1} ,\, q_j )  + {{\partial L_N}\over{\partial q_\ell}}  (q_j ,\, q_{j+1} ) = 0 .
\monend
Taking into account the  relation (\ref{lagrangien-discret}), we obtain
\moneq   \label{schema-ordre-deux} 
 {{q_{j+1} - 2\, q_j + q_{j-1}}\over{h^2}} + {1\over{2 \, m}} \, \Big[  V'_{j+1/2} + V'_{j-1/2}  \Big] = 0 
\monend
with $ \, V'_{j+1/2} \equiv   {{\partial V}\over{\partial q}} \big( {{q_j + q_{j+1}}\over2} \big) $. 
Observe first that the scheme (\ref{schema-ordre-deux}) is implicit.
Secondly, this numerical scheme is clearly consistent with the second order differential
equation (\ref{oscillateur-nonlineaire})
associated with the Lagrangian
proposed in (\ref{Lagrangien-continu}).
The momentum $ \, p_r \, $ is defined by
\moneq \label{moment} 
p_r = {{\partial L_N}\over{\partial q_r}}  (q_\ell  ,\, q_r ) .
\monend
We have $ \, p_{j+1} = m \, {{q_{j+1} - q_j}\over{h}} - {{h}\over{2}} \, V'_{j+1/2} \, $
and an analogous relation for $ \, p_j $. Then after two lines of algebra, we obtain
the two relations
$ \, p_{j+1} -  p_j = -h \,  V'_{j+1/2} \, $ and $ \, p_{j+1} +  p_j = {{2 \, m}\over{h}} \, (q_{j+1} -  q_j ) $
and a discrete system involving the momentum and the state:
\moneq \label{systeme-dynamique-discret-ordre-2} 
p_{j+1} - p_j + h \,  {{\partial V}\over{\partial q}} \Big( {{q_j + q_{j+1}}\over2} \Big) = 0 \,,\,\,
q_{j+1} - q_j -  {h\over{2\, m}} \,  \big( p_{j+1} + p_j \big) = 0 .
\monend
These relations are consistent with the first order Hamilton equations 
(\ref{hamilton-continu}).

\smallskip \noindent
The implicit scheme (\ref{systeme-dynamique-discret-ordre-2}) is implemented as follows.
We write the equations that have to be solved at each time step 
under the form $ \, F_N(p_{j+1} ,\, q_{j+1}) = 0 $. Then we consider the  Newton algorithm
$ \, (p,\, q) \longrightarrow (p^* ,\, q^* ) $:  
\moneqstar 
F_N(p ,\, q) + \dd  F_N(p ,\, q) . (p^* - p,\, q^* - q) = 0  
\monendstar 
with the initialization $ \, p = p_j ,\, q = q_j $. Observe that the inverse of the Jacobian matrix can be evaluated easily:
\moneqstar
\big( \dd F_N(p ,\, q) \big)^{-1} = {{1}\over{1 + {{h^2}\over{4\,m}} \, V'' ({{q_j+q}\over2}) }} \,
\begin{pmatrix} 1 & -{h\over2} \, V'' ({{q_j+q}\over2}) \\ {{h}\over{2\,m}} & 1  \end{pmatrix} .  
\monendstar 
In our experiments, this algorithm is converging to machine accuracy in five iterations, which is typical of Newton's algorithm.

\smallskip \noindent
By definition (see \cite{SS92}), symplectic schemes are area-preserving. Therefore,
it is sufficient to establish that the jacobian
$ \, J \equiv {{\partial p_{j+1}}\over{\partial p_j}} \, {{\partial q_{j+1}}\over{\partial q_j}} 
- {{\partial p_{j+1}}\over{\partial q_j}} \, {{\partial q_{j+1}}\over{\partial p_j}} \,$ is equal to 1.
From (\ref{systeme-dynamique-discret-ordre-2}),
with $ \,  V''_{j+1/2} \equiv  {{\partial V}\over{\partial q}} (q_{j+1/2}) \, $
and $ \, q_{j+1/2} = {1\over2} \, (q_j+q_{j+1}) $, we have
\vskip -.3 cm 
\moneqstar   \begin{array} {rcl}
\delta p_{j+1} + {{h}\over2} \,  V''_{j+1/2}  \, \delta q_{j+1} &=&
\delta p_j - {{h}\over2} \,  V''_{j+1/2}  \, \delta q_{j}   \\
-{{h}\over{2\, m}} \,  \delta p_{j+1} + \, \delta q_{j+1} &=&
{{h}\over{2\, m}} \,  \delta p_{j} + \delta q_{j} .
\end{array}  \monendstar

\smallskip \noindent 
This system can be rewritten
\moneq \label{newmark-symplectic} 
A \, \begin{pmatrix} \delta p \\ \delta q  \end{pmatrix}_{j+1} \!\! = B \, \begin{pmatrix} \delta p \\ \delta q  \end{pmatrix}_{j}  .
\monend
We have $ \, {\rm det} \, A = 1 +  {{h}^2\over{4 \, m}} \,  V''_{j+1/2}  = {\rm det} \, B \, $ and in consequence 
$ \, J = 1  $.
Observe that although the Newmark scheme  (\ref{systeme-dynamique-discret-ordre-2})  is symplectic, it does not preserve the energy.
While the quantity $ \, H_j \equiv {1\over{2\,m}} \, p_j^2 + V (q_j) \, $  remained constant for the linear version of the scheme
(see \cite{Du23}), it is not strictly constant in this nonlinear case.

\bigskip    \noindent {\bf \large    4) \quad    Simpson's quadrature  with quadratic interpolation} 

\smallskip \noindent
Internal interpolation between $ \, 0 \, $ and $ \, h \, $ is written
in terms of quadratic finite elements (see {\it e.g.}  \cite{RT83}). For $ \, 0 \leq \theta \leq 1 $, we first set 
\moneq \label{fonctions-de-base} 
\varphi_0(\theta) =  (1 - \theta) \, (1 - 2 \, \theta) \,,\,\,  
\varphi_{1/2}(\theta) =  4 \, \theta  \, (1 - \theta) \,,\,\,   
\varphi_{1}(\theta) =  \theta \, (2 \, \theta  - 1 ) .   
\monend
With $ \, t = h \, \theta $, we consider the polynomial function
\moneq \label{q2t} 
q(t) = q_\ell \,\varphi_0(\theta) + q_m \, \varphi_{1/2}(\theta) +  q_r \, \varphi_{1}(\theta) . 
\monend 
Then $ \, q(0) = q_\ell $, $ \, q({{h}\over2}) =  q_m \, $ and $ \, q(h) = q_r \, $
and the basis functions (\ref{fonctions-de-base}) are well adapted to these degrees of freedom. 
We have also $\, {{\dd q}\over{\dd t}} =  g_\ell \, (1-\theta) + g_r \, \theta \, $
with the derivatives $ \, g_\ell \, $ and $ \, g_r \, $ given by the relations 
\moneq \label{derivees-gear} 
g_\ell =  {{\dd q}\over{\dd t}}(0) =  {{1}\over{h}} \,  \big( -3 \, q_\ell + 4 \,  q_m - q_r \big) \,,\,\, 
g_r =  {{\dd q}\over{\dd t}}(h) = {{1}\over{h}} \,  \big(  q_\ell - 4 \,  q_m + 3 \, q_r \big) .
\monend 
We remark also that 
\moneq \label{derivee-centree}
 g_m = {{\dd q}\over{\dd t}} \Big( {{h}\over2} \Big) = {1\over2} \, (g_\ell + g_r) = {{q_r - q_\ell}\over{h}} .
\monend 

\smallskip \noindent 
Once the interpolation is defined in an interval of length $ \, h $, we use it by splitting
the range~$ \, [0,\,T] \, $  into $ \, N \, $  pieces, and $ \, h = {{T}\over{N}} $.
With $ \, t_j = j \, h $, we set $ \, q_j \simeq q(t_j) \, $ for $ \, 0 \leq j \leq N \, $
and $ \, q_{j+1/2} \simeq q(t_j+{{h}\over{2}}) \, $ with $ \, 0 \leq j \leq N-1 $.
In the interval $ \, [t_j ,\, t_{j+1} ] $, the function $ \, q(t) \, $ is a polynomial of degree 2, represented by the
relation (\ref{q2t}) with $ \, t = t_j + \theta \, h  $, $ \, q_\ell = q_j $, $\, q_m = q_{j+1/2} \, $
and $ \, q_{r} =  q_{j+1} $. 

\smallskip \noindent
For the numerical integration of a regular function $ \, \psi \, $ on the interval $\, [0,\, 1] $,
the Simpson method is very popular: 
$ \, \int_0^1 \psi(\theta) \, \dd \theta \simeq {1\over6} \, \big[ \psi(0) + 4 \, \psi\big( {1\over2} \big) + \psi(1) \big] $. 
%
This quadrature formula  
is exact  up to polynomials of  degree three. 
Then a discrete  Lagrangian
$ \, L_S (q_\ell ,\, q_m ,\, q_r) \simeq \int_0^h  \big[ {m\over2}  \big( {{\dd q}\over{\dd t}} \big)^2 - V(q) \big] \, \dd t \, $ 
can be defined  with the Simpson quadrature formula 
associated with an
internal polynomial approximation~$ \, q(t) \, $ of degree~2 presented in~(\ref{q2t}): 
\moneq \label{LLhh} 
L_S (q_\ell ,\, q_m ,\, q_r) =  {{m\,h}\over{12}} \, \big( g_\ell^2 + 4 \, g_m^2 + g_r^2 \big)
-  {{h}\over{6}} \, \big( V_\ell + 4 \, V_m + V_r \big)  
\monend 
with $ \, V_\ell = V(q_\ell) $, $ \, V_m = V(q_m)  \, $ and $ \, V_r = V(q_r) $. 
The discrete action $ \, \Sigma_d \, $ for a motion $ \, t \longmapsto q(t) \, $
between the initial time and a given time $ \, T > 0 \, $ is discretized
with $ \, N \, $ regular intervals and take the form 
$ \, \Sigma_d = \sum_{j=1}^{N-1} L_S (q_j ,\, q_{j+1/2} ,\, q_{j+1} ) $. 
%

\smallskip \noindent
Euler-Lagrange equations, coming from the stationary action $ \, \delta \Sigma_d = 0 \, $,
are first established for an arbitrary variation $ \, \delta  q_{j+1/2} \, $ of the internal degree of freedom
in the interval $ \, [t_j ,\, t_{j+1} ] $:
%
\moneq \label{equation-milieu-intervalle-1}
{{\partial L_S}\over{\partial q_m}} = 0 . 
\monend
Due to the relations (\ref{derivees-gear})(\ref{derivee-centree}), 
we first observe that
$ \,  {{\partial g_\ell}\over{\partial q_m}} = {{4}\over{h}} $, 
$ \,  {{\partial g_m}\over{\partial q_m}} = 0 \, $ and
$ \,  {{\partial g_r}\over{\partial q_m}} = -{{4}\over{h}} $. 
Then, due to the expression (\ref{LLhh}) of the discrete Langrangian, 
we have, with $ \, V'_m \equiv {{\partial V}\over{\partial q}} (q_m) $, 

\smallskip \noindent
$ {{\partial L_S}\over{\partial q_m}} = {{m \, h}\over6} \, \big( g_\ell \, {4\over{h}} + g_\ell \, {4\over{h}} \big) 
- {2\over3} \, h \, V'_m $
$ =  {2\over3} \, {m\over{h}} \, \big( -4 \, q_\ell + 8 \, q_m - 4 \, q_r \big) -  {2\over3} \, h \, V'_m $

\smallskip \noindent \qquad 
$ \! = {16\over3} \, {m\over{h}} \, \big[ q_m - {1\over2} \, (q_\ell + q_r) - {{h^2}\over{8 \, m}} \, V'_m \big] \,\, $ 
and the equation (\ref{equation-milieu-intervalle-1}) can be written
\moneq \label{equation-milieu-intervalle}
q_m - {{h^2}\over{8 \, m}} \,  {{\partial V}\over{\partial q}} (q_m) = {1\over2} \, (q_\ell + q_r) . 
\monend 
It defines implicitly the value $ \, q_m \, $ at the middle of the interval
as a function of the extremities~$ \, q_\ell \, $ and $ \, q_r $.
Under the form
$ \, {{4 \, m}\over{h^2}} \, (q_\ell - 2 \, q_m + q_r ) +  {{\partial V}\over{\partial q}} (q_m) = 0 $,
the relation (\ref{equation-milieu-intervalle}) is clearly consistent with the 
differential equation (\ref{oscillateur-nonlineaire}). 

\smallskip \noindent
The discrete action takes the form
$ \, \Sigma_d = \dots + L_S (q_{j-1}  ,\, q_{j-1/2} ,\, q_j) + L_S (q_{j}  ,\, q_{j+1/2} ,\, q_{j+1}) + \dots $.
The variation of this discrete action is equal to zero. We have in consequence
the following discrete Euler-Lagrange equations: 
$ \, {{\partial L_S}\over{\partial q_r}}  (q_{j-1}  ,\, q_{j-1/2} ,\, q_j)   + 
 {{\partial L_S}\over{\partial q_\ell}} ( q_{j}  ,\, q_{j+1/2} ,\, q_{j+1} ) = 0 $. 
%
We obtain after some lines of elementary calculus
\moneq \label{euler-lagrange-sommet}
{{1}\over{h^2}} \, \big( q_{j-1} -2\, q_j + q_{j+1} \big) + {1\over{3\, m}} \, \big( V'_{j-1/2} + V'_j +  V'_{j+1/2} \big) = 0  .
\monend

\smallskip \noindent
Due to the condition (\ref{equation-milieu-intervalle-1}), the right momentum $ \, p_r = {{\partial L_S}\over{\partial q_r}} \, $
can be evaluated as follows:

\smallskip \noindent
$  p_r = {{m \, h}\over6} \, \big[ g_\ell \, \big(-{{1}\over{h}} \big) + 4 \, g_m \, \big({{1}\over{h}} \big) + g_r \big({{3}\over{h}} \big) \big] - {h\over6} \, V'_r   $ 

\smallskip \noindent \quad 
$ \, = {{m}\over6} \,  \big[ -{{1}\over{h}} \, \big( - 3 \, q_\ell + 4 \, q_m - q_r \big) +  {{1}\over{h}} \, (q_r - q_\ell )
+ {{3}\over{h}} \, \big( q_\ell - 4 \, q_m + 3 \, q_r \big) \big] - {h\over6} \, V'_r   $
  
\smallskip \noindent \quad 
$ \, =  {{m}\over{6\,h}} \,  \big[ 2 \, q_\ell - 16 \, q_m + 14 \, q_r \big] - {h\over6} \, V'_r   $
$ =  {{1}\over{6}} \,  \big[ {{2 \, m}\over{h}} \, \big( q_\ell - 8 \, q_m + 7 \, q_r \big) - h  \, V'_r \big]   $.

\smallskip \noindent
We replace the intermediate value $ \, q_m \, $ with the relation (\ref{equation-milieu-intervalle}) and we have 

\smallskip \noindent
$  p_r =  {{1}\over{6}} \,  \big[ {{2 \, m}\over{h}} \, \big( q_\ell
- 8 \, \big( {1\over2} \, (q_\ell + q_r) + {{h^2}\over{8 \, m}} \,  V'_m \big)  + 7 \, q_r \big) - h  \, V'_r \big] $
$ =  {{m}\over{h}} \, \big(q_r - q_\ell \big) -  {{h}\over{6}} \, \big( 2 \,  V'_m + V'_r \big) $.

\smallskip \noindent
In other words, 
\moneq \label{p_jp1}
p_{j+1} =   {{m}\over{h}} \, \big(q_{j+1} - q_j \big) -  {{h}\over{6}} \, \big( 2 \, V'_{j+1/2} + V'_{j+1}  \big) .
\monend  
Similarly, we have
$ \, p_{j} =   {{m}\over{h}} \, \big(q_{j} - q_{j-1} \big) -  {{h}\over{6}} \, \big( 2 \, V'_{j-1/2} + V'_{j}  \big) \, $
and taking into account the discrete Euler-Lagrange equations (\ref{euler-lagrange-sommet}), we find 
\moneq \label{p_j}
p_{j} =   {{m}\over{h}} \, \big(q_{j+1} - q_j \big) +  {{h}\over{6}} \, \big( V'_j + 2 \, V'_{j+1/2} \big) .
\monend  
From the relations (\ref{p_jp1}) and (\ref{p_j}), we deduce the discrete Hamiltonian dynamics
\moneq \label{simpson-hamilton} \left\{  \begin{array} {l}
p_{j+1} - p_j +  {{h}\over{6}} \, \big( V'_j + 4 \,  V'_{j+1/2} + V'_{j+1}  \big) = 0 \\
q_{j+1} - q_j -  {{h^2}\over{12 \, m}} \, \big(  V'_{j+1} -  V'_{j} \big) - {{h}\over{2 \, m}} \, \big( p_{j+1} + p_j \big) = 0 
\end{array} \right. \monend  
We write the system (\ref{equation-milieu-intervalle})(\ref{simpson-hamilton}) under the form
$ \, F_S (q_{j+1/2} ,\, p_{j+1} ,\, q_{j+1} ) = 0 $. We have
\moneqstar
\dd F_S (q_m ,\, p ,\, q ) = \begin{pmatrix} 1 - {{h \, \theta}\over8} & 0 & - {1\over2} \\
{2\over3} \, m \, \theta & 1 & {1\over6} \, m \, \varphi \\ 0 &  -{{h}\over{2 \, m}} & 1 -  {{h \, \varphi}\over12} \end{pmatrix} 
\monendstar 
with $ \, \theta \equiv  {{h}\over{m}} \, V''(q_m) \, $ and $ \, \varphi \equiv   {{h}\over{m}} \, V''(q) $. 
After a  formal calculation with the help of the free software ``SageMath'' \cite{sagemath}, we explicit the inverse of this
Jacobian matrix:
\moneqstar \left\{  \begin{array} {l}
\big( \dd F_S (q_m ,\, p ,\, q ) \big)^{-1} \\ = \begin{pmatrix}
1 & {{1}\over{4}} \, {{h}\over{m}} & {1\over2} \\ -{2\over3} \, m \, \theta + {1\over18} \, h  \, m  \, \theta  \, \varphi &
1 - {1\over8} \, h  \,  \theta  - {1\over12} \, h \, \varphi +  {1\over96} \, h^2  \, \theta  \, \varphi &
-{1\over3} \, m \, \theta -{1\over6} \, m \, \varphi +   {1\over48} \, h  \, m  \, \theta  \, \varphi \\
-{1\over3} \, h \, \theta & {1\over2} \, {{h}\over{m}} - {1\over16} \, {{h^2 \, \theta}\over{m}} & 1 -  {1\over8} \, h \, \theta
\end{pmatrix} . 
\end{array} \right. \monendstar 
%
At each time step, the numerical resolution of the nonlinear system of three equations
$ \,  F_S (q_{j+1/2} ,\, p_{j+1} ,\, q_{j+1} ) = 0 \, $ is conducted with a Newton algorithm.
As with Newmark's algorithm, we have observed machine precision convergence at the fifth
iteration with the proposed scheme.

\smallskip \noindent
The Simpson scheme  (\ref{equation-milieu-intervalle})(\ref{simpson-hamilton})
is symplectic. 
With the notations $ \,  V''_{j} \equiv  {{\partial^2 V}\over{\partial q^2}} (q_j) \, $
and\br
$ \, \widetilde{V}''_{j+1/2} \equiv {{2 \, V''_{j+1/2}}\over{1 - {{h}^2\over{8 \, m}} \,  V''_{j+1/2}}}  $, 
we have by differentiation of  the relations (\ref{simpson-hamilton}):
\vskip -.1 cm 
\moneqstar   \begin{array} {rcl}
\delta p_{j+1} + {{h}\over6} \, ( V''_{j+1} + \widetilde{V}''_{j+1/2} ) \, \delta q_{j+1} &=&
\delta p_j - {{h}\over6} \, ( V''_{j} + \widetilde{V}''_{j+1/2}  ) \, \delta q_{j}   \\
-{{h}\over{2\, m}} \,  \delta p_{j+1} + \big( 1 -  {{h}^2\over{12 \, m}} \,  V''_{j+1}  \big)  \, \delta q_{j+1} &=&
{{h}\over{2\, m}} \,  \delta p_{j} + \big( 1 -  {{h}^2\over{12 \, m}} \,  V''_{j}  \big)  \, \delta q_{j} .
\end{array}  \monendstar
As for the Newmark scheme, this system can be written with a relation (\ref{newmark-symplectic}). We have  in this case
$ \, {\rm det} \, A = 1 +  {{h}^2\over{12 \, m}} \, \widetilde{V}''_{j+1/2} = {\rm det} \, B \, $ 
and the relation 
$ \, {{\partial p_{j+1}}\over{\partial p_j}} \, {{\partial q_{j+1}}\over{\partial q_j}} 
- {{\partial p_{j+1}}\over{\partial q_j}} \, {{\partial q_{j+1}}\over{\partial p_j}} = 1 \, $
is established.


\bigskip  \noindent {\bf \large 5) \quad   First numerical experiments and conclusions} 

%
\smallskip \noindent
We have implemented the Simpson symplectic scheme  (\ref{equation-milieu-intervalle})(\ref{simpson-hamilton}) and
have compared it with the Newmark scheme (\ref{systeme-dynamique-discret-ordre-2}). 
Typical results for $ \, N=10 \, $ meshes and one period are displayed in Figure~2.
They are compared with the exact solution presented in Figure 1. 
Quantitative errors with the maximum norm are also presented in Table~1 below.
An asymptotic order of convergence can be estimated for the momentum, the state and various energies.

\bigskip
\centerline {\includegraphics[width=1.0\textwidth]{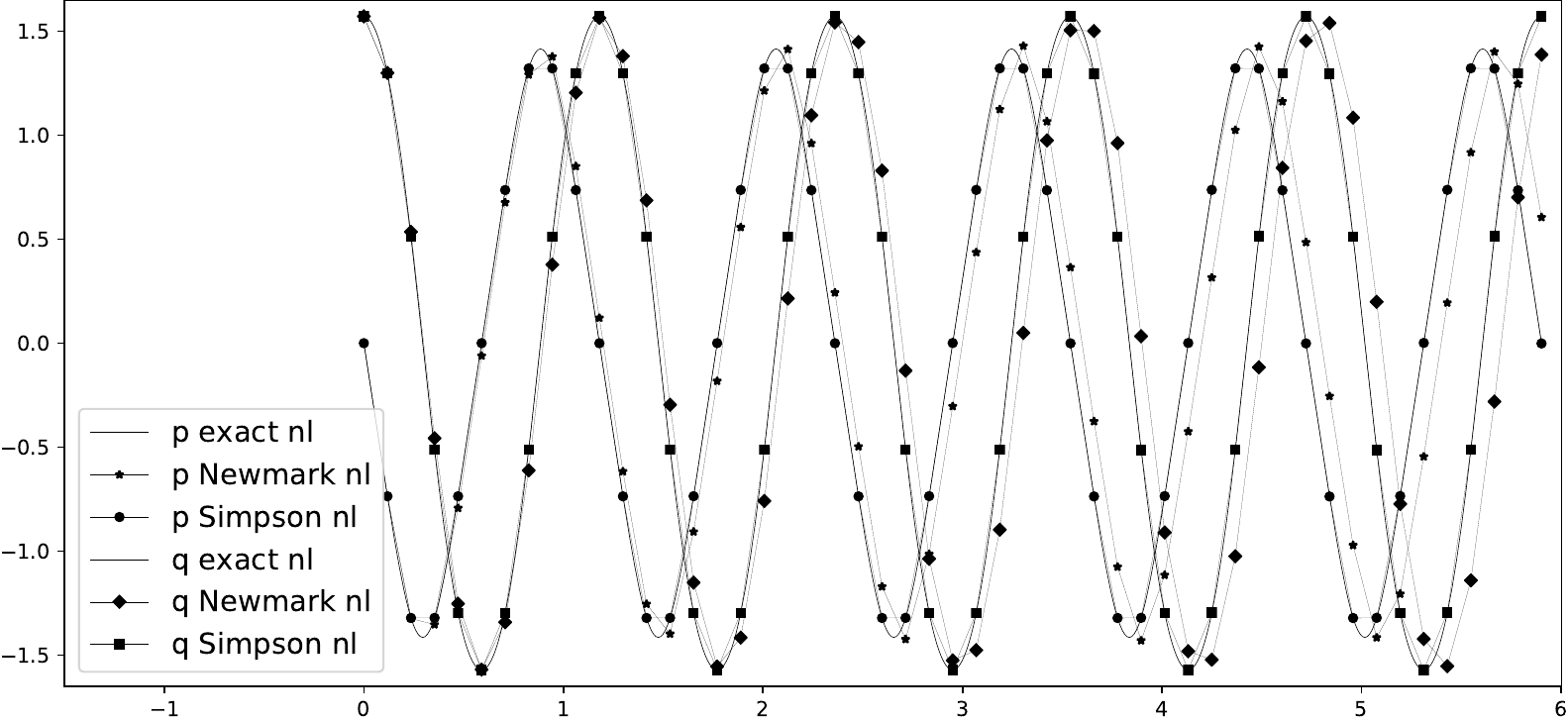}}

\smallskip  \noindent
Figure 2. Comparing Newmark and Simpson schemes with  10 points per period and a total time of 5 periods.
The exact solution is presented in Figure 1. 
The errors for the  Newmark scheme are clearly visible whereas 
the Simpson scheme is still very precise despite the reduced number of time steps. 

\bigskip 
\centerline { \begin{tabular}{|c|c|c|c|c|c|c|c|c|c|c|}    \hline
& number of meshes & 50 & 100 & 200 & order \\   \hline
Newmark & momentum & $ 2.93 \,\, 10^{-2} $  &  $ 7.32 \,\, 10^{-3} $  & $ 1.83 \,\, 10^{-3} $  & 2.0  \\   \hline
Symplectic Simpson & momentum & $  6.08 \,\, 10^{-6} $  &  $  3.78 \,\, 10^{-7} $   &  $   2.36  \,\, 10^{-8} $        & 4.0 \\   \hline
Newmark & state    & $ 5.26  \,\, 10^{-3} $  &  $ 1.31  \,\, 10^{-3} $  & $ 3.29   \,\, 10^{-4} $  & 2.0   \\   \hline
Symplectic Simpson  & state    & $ 1.05   \,\, 10^{-6} $  & $ 6.51  \,\, 10^{-8} $   & $ 4.06  \,\, 10^{-9} $     & 4.0   \\   \hline
Newmark & relative energy &  $ 9.06 \,\, 10^{-4} $ &   $ 2.29 \,\, 10^{-4} $  & $ 5.73  \,\, 10^{-5} $   & 2.0  \\   \hline
Symplectic Simpson  & relative energy  &  $ 1.30  \,\, 10^{-6}  $ &  $  8.42 \,\, 10^{-8}  $ &  $ 5.25 \,\, 10^{-9}  $ & 4.0   
\\   \hline \end{tabular} }

\smallskip  \smallskip \smallskip \noindent
Table 1. Errors in the maximum norm {($L_\infty$) for a simultation of a single period over all time steps.
For a given discretization level, the Simpson scheme is more precise than the Newmark scheme by approximately three orders of magnitude.
The orders of convergence are preliminary estimates.

\smallskip  \bigskip \noindent
In this work we have recalled Newmark's classical method which is very popular in some fields of the engineering sciences
and is symplectic. We have proposed an alternative symplectic variational integrator based on Simpson's rule,
to deal with nonlinear differential equations. The method was tested on the non-trivial nonlinear pendulum
for which the analytical solution is known. We will be working on more complex problems in the future, starting
with the symmetric spinning top.
The authors thank the reviewers for their valuable comments.

\bigskip   \newpage   \noindent {\bf \large     Annex. Exact solution of the nonlinear pendulum} 

\smallskip \noindent
In this Annex, we follow essentially the synthesis \cite{Ch2000}.
We consider an angle $ \, \theta_0 \in (0,\, \pi) \, $ and the non linear pendulum problem
\moneq \label{edo-pendule} 
  {{\dd^2 q}\over{\dd t^2}} + \omega^2 \, \sin q = 0 ,\,\,   
  q(0) = \theta_0  , \,\,  {{\dd q}\over{\dd t}}(0)  = 0 .   
\monend 
It is easy to verify the conservation of energy:
$ \, {{\dd}\over{\dd t}} \big[ {1\over2} \big(  {{\dd q}\over{\dd t}} \big)^2
   + \omega^2 \, (1 - \cos q) \big] = 0 $.
%
We introduce the parameter
$ \, k \equiv \sin \, \big( {{\theta_0}\over{2}} \big) \, $ 
%
and we have $ \, 0 < k < 1 $. Due to the initial conditions in (\ref{edo-pendule}), we observe that
$ \, {{\dd q}\over{\dd t}} < 0 \, $ for small values of $ \, t > 0 $. 
Then  we have 
$ \, {{\dd q}\over{\dd t}}  = - 2 \, \omega \, k \, \sqrt{1 - {{1}\over{k^2}} \, \sin^2 \big( {{q}\over{2}} \big) } $. 
%
We introduce also the change of unknown  defined by
$ \, \sin \varphi = {{1}\over{k}} \, \sin \big( {{q}\over2} \big) $. 
%
We observe that $ \, \varphi(0) = {{\pi}\over2} $. 
We have also
$ \,  {{\dd \varphi}\over{\dd t}} $ 
$ = -\omega \, \sqrt{1 - k^2 \, \sin^2 \varphi} = - \omega \, \cos \big( {{q}\over2} \big) \, $ 
and we deduce the differential relation
\moneq  \label{d-temps}
\omega \, \dd t = -{{\dd \varphi}\over{\sqrt{1 - k^2 \,  \sin^2 \varphi}}} 
\monend
from the previous calculus.

\smallskip \noindent 
We introduce the incomplete elliptic integral  of the first kind
$ \, F(a,\,m)  \equiv \int_0^a {{\dd \varphi}\over{\sqrt{1- m\,\sin^2 \varphi}}} \, $
and the complete elliptic integral of the first kind
$ \, K(m) \equiv \int_0^{\pi/2} {{\dd \varphi}\over{\sqrt{1- m\,\sin^2 \varphi}}} = F({{\pi}\over2},\,m)  \, $
presented {\it e.g.} in the book of  Abramowitz and Stegun \cite{AS64}.
Then by integration of  (\ref{d-temps}), 
\moneq \label{temps}
\omega \, t = - \int_{\pi/2} ^\varphi {{\dd \xi}\over{\sqrt{1 - k^2 \,  \sin^2 \xi}}} = K(k^2) - F(\varphi,\, k^2) . 
\monend
After  a fourth of a period $\, T $, the parameter $ \, \varphi \, $ changes from  $ \, {{\pi}\over2} \, $
to zero and we obtain 
\moneq  \label{periode}
\omega \, {{T}\over4} = \int_0^{\pi/2} {{\dd \varphi}\over{\sqrt{1- k^2 \,\sin^2 \varphi}}} = K(k^2) .
\monend
The Jacobi amplitude  $ \, A( . ,\, m) \, $ is defined (see  {\it e.g.} \cite{AS64}) as 
the reciprocal function of the  incomplete elliptic integral  of the first kind: 
$ \,  F(a,\,m) = u \, $ 
is equivalent to  $\, a = A (u,\, m) $. 
Due to~(\ref{temps}),   we have
$\, \varphi = A ( K(k^2) - \omega \, t ,\, k^2) $; 
finally,  $ \, \sin \big( {{q}\over2} \big) = k \, \sin \varphi \, $ and
$ \,  {{\dd q}\over{\dd t}} =  - 2 \, \omega \, k \, \cos \varphi $.

\bigskip  \bigskip     \noindent {\bf  \large  References }


\end{document}